\documentclass[12pt]{article}
\usepackage{amssymb}
\usepackage{epsfig}

\newtheorem{theorem}{Theorem}[section]
\newtheorem{lemma}[theorem]{Lemma}
\newtheorem{coro}[theorem]{Corollary}
\newtheorem{proposition}[theorem]{Proposition}

\newtheorem{definition}[theorem]{Definition}

\newtheorem{assumption}{Assumption}
\newcommand{\Z}{\mathbb Z}

\newcommand{\om}{\omega}

\newcommand{\R}{\mathbb R}
\newcommand{\C}{\mathbb C}
\newcommand{\N}{\mathbb N}

\def\tr{\rm tr}
\def\X{\cal X}
\def\ov{\overline}

\begin{document}
\title{Quantum multiple construction of subfactors}
\author{Marta Asaeda\footnote{partially supported by the NSF grant \#DMS-0202613 }}
 \date{}
 \maketitle
\begin{abstract}
 We construct the quantum $s$-tuple subfactors for an AFD II$_{1}$ subfactor with finite index and depth, for an arbitrary natural number $s$. This is a generalization of the quantum multiple subfactors by J.Erlijman and H.Wenzl \cite{HJ}, which  generalizes  the quantum double construction of a subfactor for the case that the original subfactor gives rise to a braided tensor category. In this paper we give a multiple construction for a subfactor with a weaker condition than braidedness of the bimodule system. 
%
%
 \end{abstract}

\section{Introduction}
The asymptotic subfactors of AFD II$_{1}$ subfactors with finite index and depth was  constructed by A.Ocneanu \cite{O2} \cite{O3}, and S.Popa \cite{P17}, which is regarded as  Drinfel'd's quantum double construction in the language of subfactor theory. J.Erlijman gave a multiple construction for subfactors arising from braid group representations, which generalizes the  double construction for a certain class of subfactors. Further, she and H.Wenzl gave the multiple construction for braided categories, which includes the cases for subfactors that give rise to braided categories, and obtained the dual principal graphs for several cases (\cite{HJ}).   In this paper we construct the quantum multiple subfactors for subfactors whose paragroup satisfies the generalized Yang-Baxter equation \cite{K6}. The class of subfactors with this condition includes the ones with non-braided bimodule system such as type $E_6$, $E_8$ subfactors. It also includes subfactors with non-commutative bimodule system, such as $M \subset M \rtimes {\frak S}_3$. It is expected that the  quantum multiple subfactors constructed in this paper are of finite depth. It is easily observed that the subfactors constructed in this paper include the ones given in \cite{HJ}.

Throughout this paper all the von Neumann algebras are of type AFD, and all the subfactors are assumed to be of finite index and finite depth, except the ones that we are about to construct, for which these properties need to be proved.  For the definitions of the terms such as paragroups, connections, flatness, string algebras,   see \cite{EK7} Ch.9-11. 
\section{The commuting square and biunitary connection}
Let $N \subset M$ be a subfactor and $(G, H, \beta, W)$ be its paragroup, as in the following picture, where $G$ (resp. $H$) is the (dual) principal graph, and $\beta^{2}$ is the index of the subfactor.  Note that the notation here is upside down from the usual notation. When a graph is laid so the even vertices are on the bottom (resp. left) and the odd vertices are on the top (resp. right), we consider it as in the ``right position'', and if it is the other way around, we consider it as renormalized one and call them $W_{1}$, $W_{2}$, $W_{3}$ for horizontal renormalization, vertical renormalization, and the sequence of the two, respectively. 

\begin{center}
\thinlines
\unitlength 1.0mm
\begin{picture}(30,20)(0,-7)
\multiput(11,-2)(0,10){2}{\line(1,0){8}}
\multiput(10,7)(10,0){2}{\line(0,-1){8}}
\multiput(10,-2)(10,0){2}{\circle*{1}}
\multiput(10,8)(10,0){2}{\circle*{1}}
\put(6,12){\makebox(0,0){$V_2$}}
\put(24,12){\makebox(0,0){$V_3$}}
\put(6,-6){\makebox(0,0){$V_0$}}
\put(24,-6){\makebox(0,0){$V_1$}}
\put(15,12){\makebox(0,0){$H^{t}$}}
\put(15,-6){\makebox(0,0){$G$}}
\put(6,3){\makebox(0,0){$G$}}
\put(24,3){\makebox(0,0){$H^{t}$}}
\put(15,3){\makebox(0,0){$W$}}
\end{picture} 
\end{center}
Let $\om$ be its global index. We construct the $s$-dimensional nested graphs in the first $2^{s}$-ant of $\R^{s}$ as follows. \\ \ \\
First we construct an enlarged biunitary connection obtained as a product of $W$ and its renormalizations. 
\begin{center}
{\psfig{figure=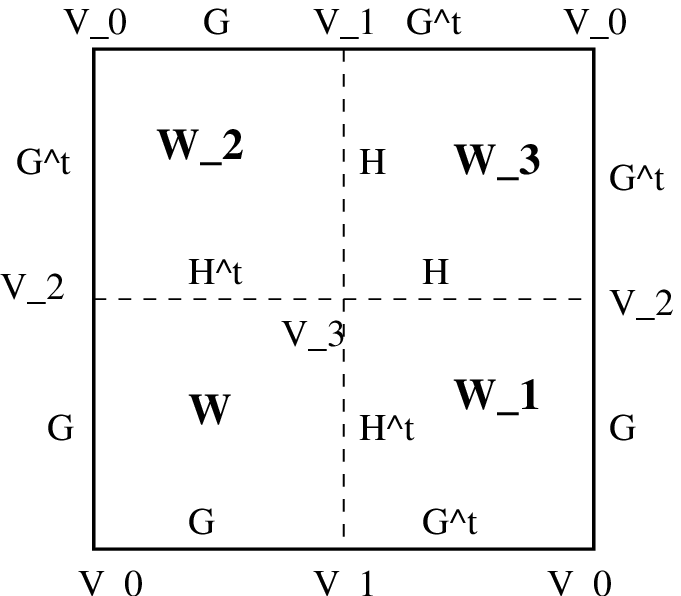,height=4cm}}
\end{center}
Let $K=G \cdot G^{t}$, i.e. a bipartite graph whose even and odd vertices are both $V:=V_{0}$, and the edges are given by concatenation of the edges in $G$ and those in $G^{t}$. We define a new biunitary connection $Y$ as follows: 
\begin{center}
\thinlines
\unitlength 1.0mm
\begin{picture}(30,20)(0,-7)
\multiput(11,-2)(0,10){2}{\line(1,0){8}}
\multiput(10,7)(10,0){2}{\line(0,-1){8}}
\multiput(10,-2)(10,0){2}{\circle*{1}}
\multiput(10,8)(10,0){2}{\circle*{1}}
\put(6,12){\makebox(0,0){$V$}}
\put(24,12){\makebox(0,0){$V$}}
\put(6,-6){\makebox(0,0){$V$}}
\put(24,-6){\makebox(0,0){$V$}}
\put(15,12){\makebox(0,0){$K$}}
\put(15,-6){\makebox(0,0){$K$}}
\put(6,3){\makebox(0,0){$K$}}
\put(24,3){\makebox(0,0){$K$}}
\put(15,3){\makebox(0,0){$Y$}}
\end{picture} 
\end{center}
\begin{eqnarray*}
&
\raisebox{-.25in}{
 \thinlines
\unitlength 1.0mm
\begin{picture}(30,20)(0,-7)
\multiput(11,-2)(0,10){2}{\line(1,0){8}}
\multiput(10,7)(10,0){2}{\line(0,-1){8}}
\multiput(10,-2)(10,0){2}{\circle*{1}}
\multiput(10,8)(10,0){2}{\circle*{1}}
\put(15,12){\makebox(0,0){$\xi_{3} \cdot \eta_{3}$}}
\put(15,-6){\makebox(0,0){$\xi_{0} \cdot \eta_{0}$}}
\put(3,3){\makebox(0,0){$\xi_{1} \cdot \eta_{1}$}}
\put(27,3){\makebox(0,0){$\xi_{2} \cdot \eta_{2}$}}
\put(15,3){\makebox(0,0){$Y$}}
\end{picture}
} 
\; :=\sum_{\nu_{i}} \raisebox{-.25in}{ {\psfig{figure=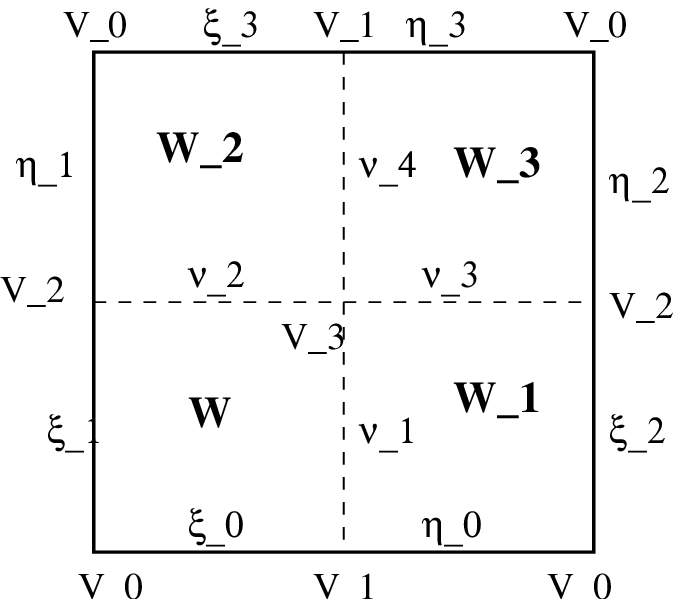,height=2.5cm}}} \\
&= \sum_{\nu_{i}}
\raisebox{-.25in}{
\thinlines
\unitlength 1.0mm
\begin{picture}(30,20)(0,-7)
\multiput(11,-2)(0,10){2}{\line(1,0){8}}
\multiput(10,7)(10,0){2}{\line(0,-1){8}}
\multiput(10,-2)(10,0){2}{\circle*{1}}
\multiput(10,8)(10,0){2}{\circle*{1}}
\put(15,12){\makebox(0,0){$\nu_{2}$}}
\put(15,-6){\makebox(0,0){$\xi_{0}$}}
\put(6,3){\makebox(0,0){$\xi_{1}$}}
\put(24,3){\makebox(0,0){$\nu_{1}$}}
\put(15,3){\makebox(0,0){$W$}}
\end{picture} 
\hspace{-1cm}
\thinlines
\unitlength 1.0mm
\begin{picture}(30,20)(0,-7)
\multiput(11,-2)(0,10){2}{\line(1,0){8}}
\multiput(10,7)(10,0){2}{\line(0,-1){8}}
\multiput(10,-2)(10,0){2}{\circle*{1}}
\multiput(10,8)(10,0){2}{\circle*{1}}
\put(15,12){\makebox(0,0){$\nu_{3}$}}
\put(15,-6){\makebox(0,0){$\eta_{0}$}}
\put(6,3){\makebox(0,0){$\nu_{1}$}}
\put(24,3){\makebox(0,0){$\xi_{2}$}}
\put(15,3){\makebox(0,0){$W_{1}$}}
\end{picture} 
\hspace{-1cm}
\thinlines
\unitlength 1.0mm
\begin{picture}(30,20)(0,-7)
\multiput(11,-2)(0,10){2}{\line(1,0){8}}
\multiput(10,7)(10,0){2}{\line(0,-1){8}}
\multiput(10,-2)(10,0){2}{\circle*{1}}
\multiput(10,8)(10,0){2}{\circle*{1}}
\put(15,12){\makebox(0,0){$\xi_{3}$}}
\put(15,-6){\makebox(0,0){$\nu_{2}$}}
\put(6,3){\makebox(0,0){$\eta_{1}$}}
\put(24,3){\makebox(0,0){$\nu_{4}$}}
\put(15,3){\makebox(0,0){$W_{2}$}}
\end{picture} 
\hspace{-1cm}
\thinlines
\unitlength 1.0mm
\begin{picture}(30,20)(0,-7)
\multiput(11,-2)(0,10){2}{\line(1,0){8}}
\multiput(10,7)(10,0){2}{\line(0,-1){8}}
\multiput(10,-2)(10,0){2}{\circle*{1}}
\multiput(10,8)(10,0){2}{\circle*{1}}
\put(15,12){\makebox(0,0){$\eta_{3}$}}
\put(15,-6){\makebox(0,0){$\nu_{3}$}}
\put(6,3){\makebox(0,0){$\nu_{4}$}}
\put(24,3){\makebox(0,0){$\eta_{2}$}}
\put(15,3){\makebox(0,0){$W_{3}$}}
\end{picture} 
}
\end{eqnarray*}
By construction $Y$ is a flat connection, and its renormalizations are identical to itself. 
This connection produces $N \subset M_{1}$, where $M_{1}$ is the basic construction of $N \subset M$ and thus the asymptotic inclusion is the same. 
\subsection{Nested algebras on higher dimensional lattice}
Now we construct the high-dimensional nested algebras using the connection $Y$. For two dimensional case, see \cite{EK7}, 11.3. 

For ${\bf n}:=(n_{0}... n_{s-1}) \in \Z_{\geq 0}^{s}$, let $v_{\bf n}$ be the lattice point located at ${\bf n}$. $v_{\bf n}$ may be sometimes simply denoted by ${\bf n}$, abusing the notation. Let $E_{{\bf n},i}$ be the lattice edge connecting ${\bf n}$ and ${\bf n}+e_i$, where $e_i$ is the unit $i$-th vector. (Note that we number the coordinates from zero.)   We define a nested graph ${\cal K}$ as follows: Let $V_{\bf n}$ be equal to $V$ as a set, located at ${\bf n}$. Let $K_{{\bf n},i}$ be identical to the graph $K$, lying along the lattice edge $E_{{\bf n},i}$. The vertices of $K_{{\bf n},i}$ are identified with $V_{\bf n} \cup V_{{\bf n}+e_i}$ in an obvious manner. We define a nested graph by ${\cal K}=\cup_{{\bf n}, i}  K_{{\bf n},i}$. A path of ${\cal K}$ is a concatenation of edges in ${\cal K}$. We call a path consisting of lattice edges  a lattice path in order to distinguish from a path of the graph. For a path $\xi$ we denote its length by $|\xi|$, the origin and the end by $s(\xi), r(\xi)$ respectively. The same notions for a lattice path are denoted similarly. For a path $\xi$, we denote by $[\xi]$ the lattice path that $\xi$ lies along. For ${\bf n}  \in \Z_{\geq 0}^{s}$ we denote $|{\bf n}|:=n_0+n_1+...+n_{s-1}$. 

Let ${\bf n}, {\bf m} \in \Z_{\geq 0}^{s}$ be so that ${\bf n}-{\bf m}\in \Z_{\geq 0}^{s}$, 
and $L$ be a lattice path with $s(L)={\bf m}$, $r(L)={\bf n}$, $|L|=|{\bf n}-{\bf m}|$. Let $p \in V_{\bf m}$, $q \in V_{\bf n}$. 
We define 
$${\rm Path}_{p,q;L}:={\rm span}_\C \{ \xi \ | \xi \in {\cal K}, s(\xi)=p, r(\xi)=q, [\xi]=L\}.$$
Note that given another lattice path $L'$ with the same condition, we have
\begin{eqnarray*}
&&{\rm Path}_{p,q;L} \cong {\rm Path}_{p,q;L'} \cong \\
&& {\rm Path}_{p,q}^{|{\bf n}-{\bf m}|}K:={\rm span}_\C \{\xi \ |\ \mbox{\rm path in } K, s(\xi)=p, r(\xi)=q, |\xi|={|{\bf n}-{\bf m}|}\}.
\end{eqnarray*}
We give an isomorphism ${\rm Path}_{p,q;L} \cong {\rm Path}_{p,q;L'}$ by identifying the basis, using the flat connection $Y$ as follows: 
\begin{definition}
\label{algebra}
Let $\xi \in  {\rm Path}_{p,q;L}$, $\eta \in {\rm Path}_{p,q;L'}$. For simplicity we assume that $L \cap  L'=\{v_{\bf m}, v_{\bf n}\}$. Let $S$ be a union of  squares with the edges in $\cup_{{\bf n}, i} E_{{\bf n}, i}$, so that $\partial S=L \cup L'$. Assume that $S$ is taken so that the area is minimum.  We define a conjugate-linear form by 
$$<\xi, \eta>=\sum_{\sigma} \prod_k Y(\sigma_k),$$
where  $\sigma=\cup_k \sigma_k$ is a surface that lies along with $S$ so that $\partial \sigma= \xi \cup \eta$, $\sigma_k$'s are distinct square with edges in ${\cal K}$, $Y(\sigma_k)$ is the evaluation of the flat connection $Y$ on $\sigma_k$. Namely, $<\xi, \eta>$ is given as a state sum of $Y$ taken over all possible surfaces that lie along $S$ with the boundary $\xi \cup \eta$. 
\end{definition}
\begin{definition} 
\label{form}
If $<,>$ as above is well-defined (i.e. if it does not depend on the choice of $S$) and non-degenerate, 
we define an isomorphism of  ${\rm Path}_{p,q;L} \to  {\rm Path}_{p,q;L'}$ by 
$$\xi \to \sum_{\eta: \ {\rm path} \in {\rm Path}_{p,q;L'} } <\xi, \eta>\eta.$$
We define a path space by 
$$ {\rm Path}_{p,q}:= {\rm Path}_{p,q;L}, $$
where the space does not depend on the choice of $L$ under the given isomorphism. We define an algebra at ${\bf n}$ by 
$$A_{\bf n} :=  {\rm Path}_{*,{\bf n}} \otimes {\rm Path}_{*,{\bf n}}^*, $$
where we consider $* \in V_{\bf 0}$,  ${\rm Path}_{*,{\bf n}} := \oplus_{q \in V_{\bf n}}  {\rm Path}_{*,q}$, and that the dual space is given with respect to $<,>$.  We denote an element in $A_{\bf n}$ by $(\xi, \eta)=\xi \otimes \eta^*$. Note that $( c \xi, \eta)=c (\xi, \eta) =(\xi, {\bar c}\eta)$ for $c \in \C$. 
The $*$-algebra structure is given by $(\xi, \eta)\cdot (\xi', \eta'):= \delta_{\eta, \xi'} (\xi, \eta')$ and ${\ov{c(\xi, \eta)}}:={\bar c}(\eta,\xi)$ for $c \in \C$
\end{definition} 
As it is, the conjugate-linear form is ill-defined since it depends on the choice of $S$. And it is not obvious if $<,>$ is indeed non-degenerate. In the following we address those issues. 

 In order for well-definedness, the flat connection $Y$ needs to satisfy an additional condition: 
 \begin{assumption} (Generalized Yang-Baxter equation)
\label{gybe}
A biunitary connection is said to satisfy the generalized Yang-Baxter equation (GYBE) if the following equality is satisfied: 
\vspace{-1cm}
\thinlines
\unitlength 1mm
\begin{center}
\begin{picture}(155,45)
\put(21,21){\vector(1,1){8}}
\put(21,19){\vector(1,-1){8}}
\put(31,29){\vector(1,-1){8}}
\put(31,11){\vector(1,1){8}}
\put(31,30){\vector(1,0){18}}
\put(41,20){\vector(1,0){18}}
\put(31,10){\vector(1,0){18}}
\put(51,29){\vector(1,-1){8}}
\put(51,11){\vector(1,1){8}}
\multiput(30,30)(20,0){2}{\circle*{1}}
\multiput(20,20)(20,0){3}{\circle*{1}}
\multiput(30,10)(20,0){2}{\circle*{1}}
\put(36,20){\makebox(0,0){$a_7$}}
\put(17,20){\makebox(0,0){$a_1$}}
\put(63,20){\makebox(0,0){$a_4$}}
\put(30,33){\makebox(0,0){$a_2$}}
\put(50,33){\makebox(0,0){$a_3$}}
\put(30,7){\makebox(0,0){$a_6$}}
\put(50,7){\makebox(0,0){$a_5$}}
\put(21,25){\makebox(0,0){$\xi_1$}}
\put(22,15){\makebox(0,0){$\xi_6$}}
\put(39,25){\makebox(0,0){$\xi_7$}}
\put(38,15){\makebox(0,0){$\xi_9$}}
\put(59,25){\makebox(0,0){$\xi_3$}}
\put(59,15){\makebox(0,0){$\xi_4$}}
\put(40,33){\makebox(0,0){$\xi_2$}}
\put(50,23){\makebox(0,0){$\xi_8$}}
\put(40,7){\makebox(0,0){$\xi_5$}}
\put(10,20){\makebox(0,0){$\displaystyle\sum_{a_7,\xi_7,\xi_8,\xi_9}$}}
\put(77,18){\makebox(0,0){$=\displaystyle\sum_{a_7,\xi_7,\xi_8,\xi_9}$}}
\put(91,21){\vector(1,1){8}}
\put(91,19){\vector(1,-1){8}}
\put(101,30){\vector(1,0){18}}
\put(91,20){\vector(1,0){18}}
\put(111,21){\vector(1,1){8}}
\put(111,19){\vector(1,-1){8}}
\put(101,10){\vector(1,0){18}}
\put(121,29){\vector(1,-1){8}}
\put(121,11){\vector(1,1){8}}
\multiput(100,30)(20,0){2}{\circle*{1}}
\multiput(90,20)(20,0){3}{\circle*{1}}
\multiput(100,10)(20,0){2}{\circle*{1}}
\put(114,20){\makebox(0,0){$a_7$}}
\put(87,20){\makebox(0,0){$a_1$}}
\put(133,20){\makebox(0,0){$a_4$}}
\put(100,33){\makebox(0,0){$a_2$}}
\put(120,33){\makebox(0,0){$a_3$}}
\put(100,7){\makebox(0,0){$a_6$}}
\put(120,7){\makebox(0,0){$a_5$}}
\put(91,25){\makebox(0,0){$\xi_1$}}
\put(92,15){\makebox(0,0){$\xi_6$}}
\put(100,23){\makebox(0,0){$\xi_7$}}
\put(112,15){\makebox(0,0){$\xi_9$}}
\put(129,25){\makebox(0,0){$\xi_3$}}
\put(129,15){\makebox(0,0){$\xi_4$}}
\put(110,33){\makebox(0,0){$\xi_2$}}
\put(111,25){\makebox(0,0){$\xi_8$}}
\put(110,7){\makebox(0,0){$\xi_5$}}
\end{picture}
\end{center}
We assume that our flat connection $Y$ satisfy GYBE. 
We call the geometric move in the equation ``GYBE move''.  
\end{assumption}
The paragroups that correspond to subfactors with braided bimodule systems satisfy GYBE, using the translation of the language of  flat connection into rational conformal field theory as in \cite{X1}, sec.2,  and applying  Reidemeister move III . One needs to construct $2 \times 2$ connections with braided system for the vertices. Note that the biunitary connections for ADE Dynkin diagrams, including non-flat ones, satisfy the relation (\cite{EK7}, sec.11.9). It is also a straightforward arithmetic computation of cells to check that the flat connection for   ${\frak S}_{3}$ group subfactor with the bimodule system corresponding to the group elements  also satisfy GYBE. The data for the cells are obtained from matrix entries of representations of the group, see \cite{EK7}, sec.10.6. Note that the choice of the basis of a representation only amount to gauge equivalence of the connections. 
 
 In the following we prove that GYBE is indeed sufficient condition for well-definedness of the conjugate-linear form $<\xi, \eta>$. 
 For future use we restate the axiom of biunitarity, and define a geometric move associated to it. 
 \begin{definition}
Recall the biunitarity of a connection implies 
\begin{center}
$
\sum_{\nu_{i}}$
\raisebox{-.25in}{
\thinlines
\unitlength 1.0mm
\begin{picture}(30,20)(0,-7)
\multiput(11,-2)(0,10){2}{\line(1,0){8}}
\multiput(10,7)(10,0){2}{\line(0,-1){8}}
\multiput(10,-2)(10,0){2}{\circle*{1}}
\multiput(10,8)(10,0){2}{\circle*{1}}
\put(15,12){\makebox(0,0){$\nu_{2}$}}
\put(15,-6){\makebox(0,0){$\xi_{0}$}}
\put(6,3){\makebox(0,0){$\xi_{1}$}}
\put(24,3){\makebox(0,0){$\nu_{1}$}}
\put(15,3){\makebox(0,0){$Y$}}
\end{picture} 
\hspace{-1cm}}
\raisebox{-.25in}{
\thinlines
\unitlength 1.0mm
\begin{picture}(30,20)(0,-7)
\multiput(11,-2)(0,10){2}{\line(1,0){8}}
\multiput(10,7)(10,0){2}{\line(0,-1){8}}
\multiput(10,-2)(10,0){2}{\circle*{1}}
\multiput(10,8)(10,0){2}{\circle*{1}}
\put(15,12){\makebox(0,0){$\nu_{2}$}}
\put(15,-6){\makebox(0,0){$\xi'_0$}}
\put(6,3){\makebox(0,0){$\nu_{1}$}}
\put(24,3){\makebox(0,0){$\xi'_{1}$}}
\put(15,3){\makebox(0,0){$Y$}}
\end{picture} 
}
 $= \delta_{\xi_1, \xi'_1} \delta_{\xi_0, \xi'_0}. $
\end{center}
The corresponding geometric move is as follows:
\begin{center} {\psfig{figure=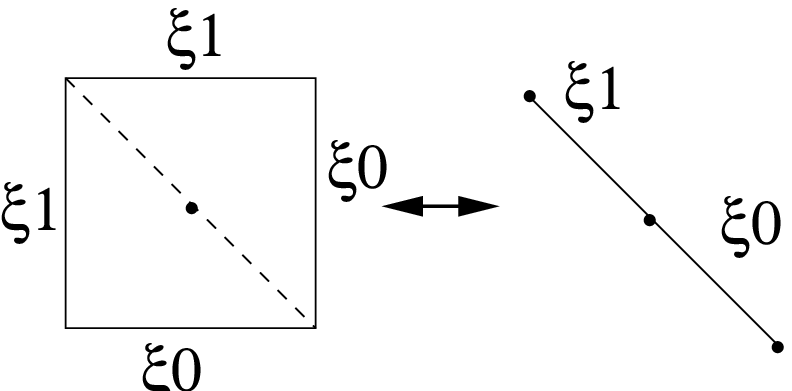,height=2.5cm}} \end{center} 
We call this unitary move. 
\end{definition}
For the time being our focus is on spacial geometry. We do not think about the actual graph ${\cal K}$ lying in the space.  We omit the word "lattice" when we discuss about paths, edges, vertices, etc. 
\begin{definition} (some notation)
\label{minsurface}
Let $v_{\bf n}$, $E_{{\bf n}, i}$ as before. We abuse notation and sometimes consider $e_{i}$ to be also an edge parallel to $e_{i}$, located possibly anywhere, i.e.  $E_{{\bf n}, i}$ for any ${\bf n}$. Consider paths $\xi$ and $\eta$ in $\cup_{{\bf n}, i} E_{{\bf n}, i}$ with $s(\xi), s(\eta)=v_{\bf n}$,  $r(\xi), r(\eta)=v_{{\bf n}+{\bf e}{\bf 1} }$, and $|\xi|=|\eta|=s$, where ${\bf e}=(e_{0}, \dots, e_{s-1})$ and ${\bf 1}=(1)_i$ (i.e. ${\bf e}{\bf 1}=e_0 + \dots + e_{s-1}$).   Notice that any such path has one-to-one correspondence with the elements of the symmetric group ${\frak S}_{s}$ by 
$$ \sigma  \in {\frak S}_{s} \leftrightarrow e_{\sigma(0)} \cdot e_{\sigma(1)} \cdot e_{\sigma(2)}  \cdot \cdots \cdot e_{\sigma(s-1)} $$
up to the initial point. 

 We consider an embedded surface in $\R^s_\geq$ to be always a union of squares with the edges in $\cup_{{\bf n}, i} E_{{\bf n}, i}$. We call a surface a minimal surface if there is no surface with a smaller area with a given boundary loop.   Note that the area is the same as the number of squares that make up the surface. 
 
 We call a surface $S$ to be ``spanned by $\xi$ and $\eta$'' if   $\partial S ={\ov{ \xi \cup \eta - \xi \cap \eta}}$ and the boundary of each component of $S$ is connected. We denote the set of the surfaces spanned by $\xi$, $\eta$ by ${\cal F}_{\xi, \eta}$, and the minimal surfaces by ${\cal MF}_{\xi, \eta}$. 
  \end{definition}
 \begin{lemma}
 \label{min}
 Any minimal surface with one boundary component  is contractible. 
 \end{lemma}
 \noindent
 {\bf Proof.} \\
  Any minimal surface with area 1 is contractible. Let $n$ be the minimum of the area for which there exists a minimal surface $S$ with a non-zero genus,  with respect to its boundary $\rho$. Let $A$ be a square in $S$ so that $A \cap \partial S$ has only one component. Such a square should exist: if all the squares with non-trivial intersection with $\partial S$ have more than one component, i.e. two opposite edges, that would give a parity among the edges in $\partial S$ that are disjoint, thus we have contradiction to the assumption that $\partial S$ is connected. 
 Consider $S \backslash A=:S'$. The area of $S'$ is smaller than that of $S$, and has the same genus. And $S'$ is a minimal surface with respect to its boundary: if there is a surface $S''$ with the same boundary with smaller area,  $S'' \cup A$ would have a smaller area than $S$, which is contradiction. Since $S'$ has a smaller area than $n$, it leads to contradiction. QED. \\ \ \\
 
 From now on we only consider the contractible surfaces and restrict the elements of  ${\cal F}_{\xi, \eta}$ to be contractible. 

\begin{proposition}
\label{gybeuni}
 Any two minimal surfaces spanned by $\xi$ and $\eta$ are deformed to each other by GYBE and unitary moves. 
 \end{proposition}
 {\bf Proof.} \\
For simplicity assume that ${\bf n}={\bf 0}$,  $\xi= e_{0} \cdot e_{1} \cdot e_{2}  \cdot \cdots \cdot e_{s-1}$. 
 and  that $\xi$ and $\eta$ do not intersect except at the beginning and the end. Let $S \in {\cal F}_{\xi, \eta}.$ $S$ is homeomorphic to a disk in this case.  We label all the edges in $S$ parallel to $e_{i}$ by $i$. Then $S$ is a union of a square looking like this:\begin{center} {\psfig{figure=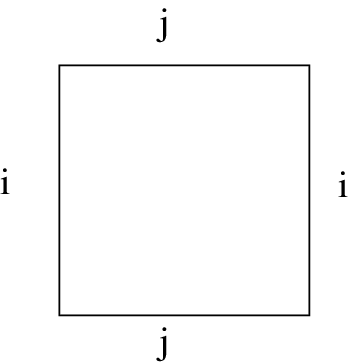,height=2.5cm}} \end{center}
 We drow a crossing on the square like this: 
 \begin{center} {\psfig{figure=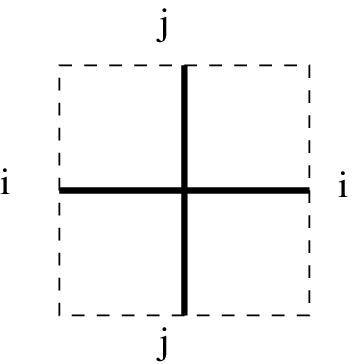,height=2.5cm}} \end{center}
 We  label the strand across $e_{i}$ also by $i$. The picture of $S$ with these decorations, for $s=4$, $\eta= e_{3} \cdot e_{2} \cdot e_{4} \cdot e_{1}$, appears as follows:
  \begin{center} {\psfig{figure=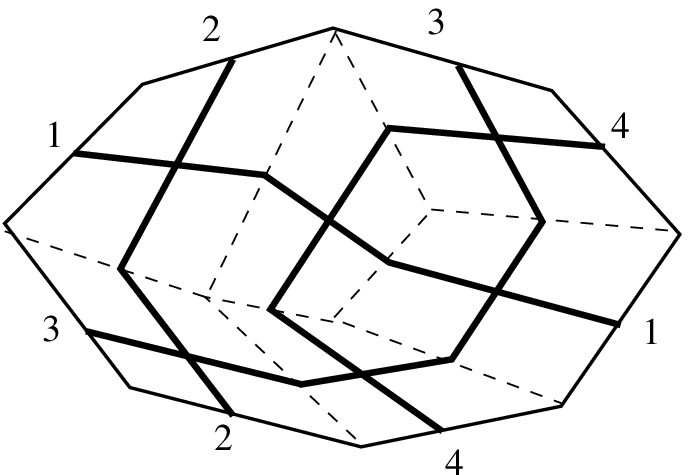,height=2.5cm}} \end{center}
  From now on we consider any surface as decorated. We observe that each $S \in {\cal F}_{\xi, \eta}$ gives rise to a degenerate tangle \\ $T_{S} \in {\cal T}_{\sigma}=$  
  \{ $ \mbox{{\rm tangles with }} s \mbox{ {\rm input and output, sending } i {\rm to }} \sigma(i)  $ without a self-crossing\}. 
 We introduce the following moves. 
\\ \ \\
{\bf Reidemeister II (R-II)} \\
\begin{center} {\psfig{figure=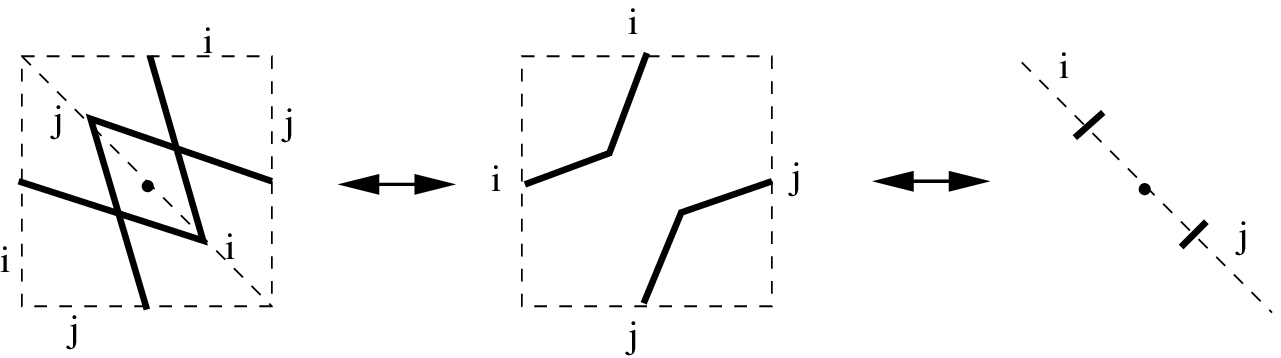,height=2.5cm}} \end{center} 
{\bf Reidemeister III (R-III)} \\
\begin{center} {\psfig{figure=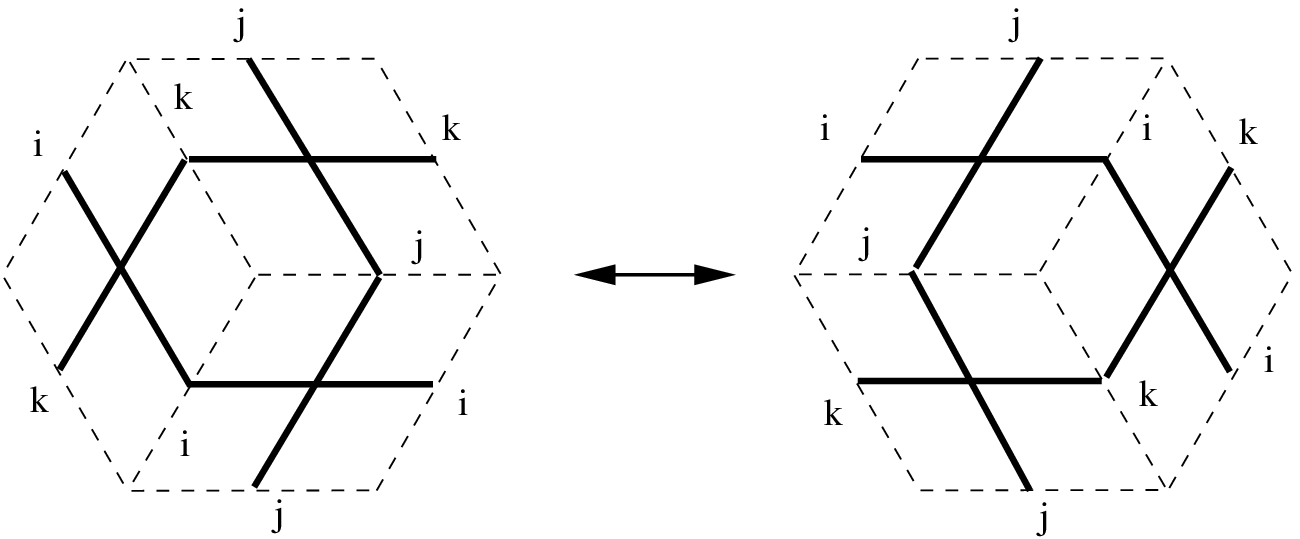,height=2.5cm}} \end{center} 
Note that those moves correspond to the Reidemeister moves in knot theory (For the definitions of Reidemiser moves as well as elementary knowledge of knot theory, see \cite{L}). The middle picture in R-II move actually does not appear in our situation but drawn just to have an association with knot theory. It immediately degenerates into the right-most picture. Note also that the Reidemeister I move does not appear in our situation since it requires a square with all four edges labeled by the same name. 
\begin{lemma}
For $T \in {\cal T}_{\sigma}$, let $c(T)$ be the number of crossings, and $c({\cal T}_{\sigma}):= {\rm min}_{T \in {\cal T}_{\sigma} } c(T)$. Then $c({\cal T}_{\sigma}) = w(\sigma)=$ the length of $\sigma$ as a word written in transpositions \{$\sigma_i=(i, i+1)$\}.  Any $T \in {\cal T}_{\sigma}$ is deformed to some $T'  \in {\cal T}_{\sigma}$ so that $c(T')=c({\cal T}_{\sigma})$ by performing R-II move and R-III moves finitely many times. 
\end{lemma}
\noindent
{\bf Proof.} The first statement is clear, as a transposition corresponds to a crossing. The second statement follows from the fact that ${\frak S}_s$ is generated by $\sigma_i$'s with relations $\sigma_i^2=1$, $\sigma_i \sigma_{i \pm 1} \sigma_i=\sigma_{i \pm 1}  \sigma_i \sigma_{i \pm 1}$, and any element $\sigma \in {\frak S}_s$ is reduced to minimum word length expression by applying those relations finitely many times.   QED. \\
%
%

The following lemma follows in a similar manner. 
\begin{lemma}
Let $T, T' \in  {\cal T}_{\sigma}$ so that $c(T)=c(T')=c( {\cal T}_{\sigma}).$ Then $T$ is deformed to $T'$ by a sequence of R-II and R-III moves. 
\end{lemma}
This lemma implies Proposition \label{gybeuni}. 

We now come back to Definition \ref{algebra}. 
For paths $\xi \in {\rm Path}_{p,q;L}$ and  $ \eta \in {\rm Path}_{p,q;L'}$, consider two minimal surfaces $S_{1}$ and $S_{2}$ spanned by $[\xi]$ and $[\eta]$. 
 Let $<\xi, \eta>_{S_{i}}$ be the conjugate-linear form given in Definition \ref{algebra}, using $S_{i}$.  
 By the above discussion, there is a sequence of intermediate surfaces $\{S_\pi\}$ that connect $S_{1}$ and $S_{2}$ where each adjacent surfaces differ by R-II and R-III moves. By the assumption of GYBE and biunitarity of the connection, $<\xi, \eta>_{S_{\pi}}$ is constant. Thus we proved that $<\xi, \eta>_{S_{1}}=<\xi, \eta>_{S_{2}}$, i.e. $<\xi, \eta>$ is well-defined. 
 
 We prove the non-degeneracy of $<\xi,\eta>$ as follows: By Lemma \ref{min}, a minimal surface $S$ spanned by $L$, $L'$ is a disk. For simplicity let us assume that $s(L)=s(L')={\bf 0}$, $r(L)=r(L')={\bf e1}$, and that $L=e_{0} \cdot e_{1}   \cdot \cdots \cdot e_{s-1}$, $L'=e_{\tau(0)} \cdot e_{\tau(1)}   \cdot \cdots \cdot e_{\tau(s-1)}$ for $\tau \in {\frak S}_s$.    The surface $S$ determines a minimal expression  $\tau=\tau_{k+1} \tau_{k}... \tau_1$, where $\tau_j$'s are transpositions. For $1 \leq j \leq k$,  Let $L_j=e_{\tau^j(0)} \cdot e_{\tau^j(1)}   \cdot \cdots \cdot e_{\tau^j(s-1)},$ where  $\tau^j=\tau_j \tau_{j-1}... \tau_1$. Let $L_0=L$, $L_{k+1}=L'$. Then we have a conjugate linear form on ${\rm Path}_{p,q, L_j} \times {\rm Path}_{p,q, L_{j+1}}$: 
$$ <\xi_j,\xi_{j+1}>= \left\{ 
\begin{array}{c}
0 \ \ \mbox{{\rm if }} \xi_j \ \mbox{and} \  \xi_{j+1} \ \mbox{{\rm disagree on } } L_j \cap L_{j+1}, \\
Y(\sigma) \  \ \mbox{{\rm otherwise,} } \phantom{XXXXXXXXXXX}
\end{array}
\right. 
$$
where $\xi_j$, $\xi_{j+1}$ are paths in ${\rm Path}_{p,q, L_j}$ and $ {\rm Path}_{p,q, L_{j+1}}$ respectively, and $\sigma$ is a square bounded by  $\xi_j$, $\xi_{j+1}$ (which corresponds to $\tau_{j+1}$). Since $Y$ is a biunitary connection, this linear form is non-degenerate. Noticing that the state sum given in Definition \ref{form} implies 
$$<\xi, \eta>=\sum_{\xi_1,... \xi_k} <\xi, \xi_1><\xi_1, \xi_2>...<\xi_k, \eta>,$$
we conclude that $<\xi,\eta>$ is non-degenerate. Note that we did not need flastness of $Y$ nor GYBE for this proof. 
 
Hence we have well-defined path spaces $ {\rm Path}_{p,q}$ and algebras $A_{\bf n}$. 
\subsection{Construction of the commuting square for the multiple subfactor}
Now we give a nested structure of algebras $\{A_{\bf n}\}$ and commuting squares arising from it. 

 We define the embedding $A_{\bf n} \subset A_{{\bf n} +e_{i}}$  by $(\xi, \eta) \mapsto \sum_{\gamma \in K_{{\bf n}, i}} (\xi \cdot \gamma, \eta \cdot \gamma)$ (i.e. $\gamma$ is parallel to $e_i$).
We define a trace on $A_{\bf n}$ by $\tr(\xi, \eta):=\delta_{\xi, \eta} \beta^{-2|{\bf n}|} \mu(r(\xi))$, where $\mu$ is the Perron-Frobenius eigenvector of the original connection $W$. This is compatible with the embedding. One can check that the conditional expectation $E_{{\bf n}, i} : A_{{\bf n}+e_{i}} \to A_{{\bf n}}$ will be given by 
$$(\xi \cdot \xi', \eta \cdot \eta') \mapsto \delta_{\xi', \eta'} \frac{\mu(r(\xi'))}{\beta^{2} \mu(r(\xi))} (\xi, \eta),$$
see \cite{EK7} Lemma 11.7. 
\begin{proposition}
\label{smallsq}
Consider the following diagram. 
$$
\begin{array}{ccc}
A_{{\bf n}+e_{j}} &\subset &  A_{{\bf n}+e_{i}+e_{j}} \\
\cup && \cup \\
A_{\bf n} &\subset &  A_{{\bf n}+e_{i}},
\end{array}
$$
where $i \neq j$. 
The identification of the bases of $A_{{\bf n}+e_{i}+e_{j}}$ via  $A_{{\bf n}+e_{i}} $ and via $A_{{\bf n}+e_{j}}$ is given by the connection $Y$. Then this is a commuting square, with conditional expectations defined as above. 
\end{proposition}
It is proved by straight forward computation: take $x=(\xi \cdot \gamma, \eta \cdot \gamma) \in A_{{\bf n}+e_{i}}$, where $(\xi, \eta) \in A_{\bf n}$. $E_{\bf n}(x)=
\frac{\mu(r(\gamma))}{\beta^{2} \mu(r(\xi))} (\xi, \eta)$. Embed $x$ into $A_{{\bf n}+e_{i}+e_{j}}$, change basis using $Y$ and apply $E_{{\bf n}+e_{j}}$; the result is  equal to $E_{\bf n}(x)$, using the unitarity of the connection. The coefficients are adjusted by the constant coming from renormalization.  QED.\\ \ \\
\begin{coro}
The following diagram is a commuting square.
$$
\begin{array}{ccc}
A_{{\bf n}+m_j e_{j}} &\subset &  A_{{\bf n}+m_i e_{i}+m_j e_{j}} \\
\cup && \cup \\
A_{\bf n} &\subset &  A_{{\bf n}+m_i e_{i}},
\end{array}
$$
where $i \neq j$, and $m_i, m_j \in \Z$. 
\end{coro}
  %
 %
\begin{theorem}
\label{multilegsq}
The following is a commuting square: 
$$
\begin{array}{ccc}
  A_{ n    {\bf e}{\bf 1}  }&\subset &  A_{(n+1)  {\bf e}{\bf 1}} \\
\cup && \cup \\
\bigvee_{i} A_{n e_{i}} &\subset &  \bigvee_{i} A_{(n+1) e_{i}},
\end{array}
$$
where the trace is defined as before, and conditional expectation is determined uniquely by the trace. 
\end{theorem}
The proof is given by successive applications of the lemma below. 
\begin{lemma}
\label{oneshift}
The following is a commuting square for all $j$: 
$$
\begin{array}{ccc}
  A_{\bf n}& \stackrel{E_{\bf n}}{\subset} &  A_{ {\bf n}+e_{j}} \\
\cup && \cup \\
\bigvee_{i} A_{n_{i} e_{i}} & \stackrel{{\tilde E_{j}}}{\subset} &  \bigvee_{i \neq j} A_{n_{i} e_{i}} \vee A_{(n_{j}+1)e_{j}},
\end{array}
$$
\end{lemma}
{\bf Proof.} \\
Using proposition \ref{smallsq} we have 
$$
\begin{array}{ccc}
  A_{\bf n}& \stackrel{E_{\bf n}}{\subset} &  A_{ {\bf n}+e_{j}} \\
\cup && \cup \\
  A_{n_{j} e_{j}} & \stackrel{{ E_{j}}}{\subset}   &   A_{(n_{j}+1)e_{j}}
\end{array} ...(**)
$$
Take $a \in A_{n_{j} e_{j}}$ and $b \in A_{(n_{j}+1)e_{j}}$. Then
 ${\tilde E_{j}}(ab)=a{\tilde E_{j}}(b)$. Thus by the uniqueness of the conditional expectation we have ${\tilde E_{j}} \mid_{A_{(n_{j}+1)e_{j}}} =E_{j}$. Now, the elements in $ \bigvee_{i \neq j } A_{n_{i} e_{i}}$ and those in $A_{(n_{j}+1)e_{j}}$  commute each other since $Y$ is a flat connection.  So it suffies to show that $ {\tilde E_{j}}(ab)  = E_{\bf n}(ab)$ for $a \in \bigvee_{i \neq j } A_{n_{i} e_{i}}$, $b \in  A_{(n_{j}+1)e_{j}}$.
 Observe that 
 \begin{eqnarray*}
  {\tilde E_{j}}(ab)&=&a{\tilde E_{j}}(b)=a E_{\bf n}(b) \ \ \mbox{\rm  by} \ (**) \\
  &=&  E_{\bf n}(ab) \ \ \mbox{\rm  because } \ a \in  A_{\bf n}.
\end{eqnarray*}
qed.   \\ \ \\
We obtain the quantum multiple inclusion $P \subset Q$ of the subfactor $N \subset M$ using the periodic commuting square as in Theorem \ref{multilegsq}, where $P:= {\ov{\cup_n  \bigvee_{i} A_{n e_{i}}}}^{w}$, $Q:=  {\ov{\cup_n   A_{ n    {\bf e}{\bf 1}  }}}^{w}$, with the inclusion given by $ \bigvee_{i} A_{n e_{i}} \subset  A_{ n    {\bf e}{\bf 1}  }$ which is compatible with the union due to the commuting square condition.   Note that when $s=2$ it coincides with the asymptotic inclusion $N \vee (N' \cap M_{\infty}) \subset M_{\infty}$.  
By comparing the commuting square given in Theorem \ref{multilegsq} and the one in \cite{HJ}, sec.3.2, one may  observe that our result here gives a generalization of \cite{HJ}. The correspondence is given by translation of the notions such as string algebras and endomorphism algebras on bimodules, embedding of an element of a string algebra and tensor product of an endomorphism with an identity map, etc. we do not discuss the detail here, however note that $X$ in \cite{HJ} is given by $_N M_N$ in this paper, and the relation between the embedding  $A_{ n    {\bf e}{\bf 1}  } \stackrel{i}{\subset} A_{ (n+1)    {\bf e}{\bf 1}  } $ in this paper and $A_{ns} \stackrel{j}{\subset} A_{(n+1)s} $ in \cite{HJ} is $i=u_{n+1} j {u_n}^*$, where $u_n$ is as defined in section 3.2 of \cite{HJ}, thus the embedding $\bigvee_{i} A_{n e_{i}}  \subset A_{ n    {\bf e}{\bf 1}  }$ in this paper does not need a unitary conjugation $u_n$ as seen in \cite{HJ}. Therefore the commuting squares in both settings are equivalent. For the correspondence between two languages, see Ch.11, 12 of \cite{EK7}. 
\section{Intermediate subfactors and the Bratteli digarams of the commuting square}
The Bratteli diagram $L$ of the inclusion $ \bigvee_{i} A_{n e_{i}} \subset A_{n {\bf e}{\bf 1}}$ is determined by the fusion structure of $N$-$N$ bimodules $\{X_{k}\}$. We show it by constructing intermediate subfactors.  
\begin{proposition}
\label{bldg}
Consider the following commuting squares: 
$$
\begin{array}{ccc}
A_{n{e}_{0}} \vee  A_{n{e}_{1}} ... \vee A_{n{e}_{s-1}} & \subset & A_{(n+1){e}_{0}} \vee  A_{(n+1){e}_{1}} ... \vee A_{(n+1){e}_{s-1}} \\
\cap && \cap \\
A_{n(e_{0}+e_{1})} \vee A_{ne_{2}} \vee ... \vee A_{n{e}_{s-1}} & \subset & A_{(n+1)({e}_{0}+e_{1})} \vee  ... \vee A_{(n+1){e}_{s-1}} \\
\cap && \cap \\
A_{n(e_{0}+e_{1}+e_{2})} \vee A_{ne_{3}} \vee ... \vee A_{n{e}_{s-1}} & \subset & A_{(n+1)({e}_{0}+e_{1}+e_{2})} \vee  ... \vee A_{(n+1){e}_{s-1}} \\
\cap && \cap \\
\vdots && \vdots \\
\cap && \cap \\
A_{n(e_{0}+ ... + e_{s-2})} \vee A_{ne_{s-1}} &\subset &  A_{(n+1)(e_{0}+ ... + e_{s-2})} \vee A_{(n+1)e_{s-1}} \\
\cap && \cap \\
A_{n {\bf 1} \cdot {\bf e}} &\subset & A_{(n+1) {\bf 1} \cdot {\bf e}}.
\end{array}
$$
The commuting square on the $j$-th floor (in European way) gives the subfactor
$$ N \vee (N' \cap M_{\infty}) \otimes \underbrace{N  \otimes ... \otimes N}_{\mbox{$j$ times}} \subset M_{\infty} \otimes \underbrace{N  \otimes ... \otimes N}_{\mbox{$j$ times}}, $$
where the embedding  is given by the asymptotic inclusion  of $N \subset M$ tensored with the identities of $N$. 
\end{proposition}
This is proved by the following lemma.
\begin{lemma}
For each $j$, the commuting square
$$
\begin{array}{ccc}
A_{n(e_{0}+ ... + e_{j-1})} \vee A_{ne_{j}} & \subset & A_{(n+1)(e_{0}+ ... + e_{j-1})} 
\vee A_{(n+1)e_{j}} \\
\cap && \cap \\
A_{n(e_{0}+ ... + e_{j})}  & \subset & A_{(n+1)(e_{0}+ ... + e_{j})} 
\end{array}
$$
gives the asymptotic inclusion. 
\end{lemma}
{\bf Proof}. \\ \ \\
Let $B_{m,n}:= A_{m(e_{0}+... + e_{j-1}), ne_{j}}$ Then the above commuting square is written as follows. 
$$
\begin{array}{ccc}
B_{n,0} \vee B_{0,n} & \subset &  B_{n,0} \vee B_{0,n} \\
\cap && \cap \\
B_{n,n}  & \subset & B_{n+1,n+1} \\
\end{array}
$$
This commuting square gives the asymptotic subfactor of $B_{0,\infty} \subset B_{1,\infty} $
Since the commuting net of algebras $ \{ B_{n,m} \}$ is given by the biunitary connection
$$\underbrace{Y Y Y ... Y}_{\mbox{$m$ times}}$$
(composed horizontally) which gives the same subfactor $N \subset M_{1}$ as $Y$ does,  we obtain the asymptotic inclusion. 
qed.
  \\
  
Thus the Bratteli diagram in each step of the left column of the diagram in Proposition \ref{bldg} is given by $\{ fusion \ graph\ \times  \times^{j} \ trivial \ graph \}$. Connecting all of this, the Bratteli diagram $L$ of the inclusion $ \bigvee_{i} A_{n e_{i}} \subset A_{n {\bf e}{\bf 1} }$ is given by the $s$-fusion graph of $N$-$N$ bimodules $\{X_{k}\}$, that is,  the set of vertices corresponding to the simple components in $\bigvee_{i} A_{n e_{i}}$ is given by $\{ (X_{0},..., X_{s-1}) \}_{X_{j} \in {\cal X}_{N-N}}$, the set of vertices corresponding to the simple components in $ A_{n{\bf 1} \cdot {\bf e}}$ is given by ${\cal X}_{N-N}$, and the number of edges between $ (X_{0},..., X_{s-1})$ and $Y$ is given by $N_{X_{0},..., X_{s-1}}^{Y}:={\rm dim Hom}(X_{0} \otimes_{N}...\otimes_{N} X_{s-1}, Y)$. In particular this implies that $P  \subset Q$ is irreducible if $N \subset M$ is irreducible.  
\\ 

The following proposition is obtained directly from the construction. 
\begin{proposition} 
The Perron-Frobenius eigenvalue $\beta_{L}$ of $L$ is given by $\om^{\frac{s-1}{2}}$, and the Perron Frobenius eigenvector $\mu_{L}$ is given by  . 
\begin{eqnarray*}
\mu_{L}(i_{0}, i_{1}, ... , i_{s-1})&=& \mu(0)\mu(1)...\mu(s-1), \\
\mu_{L}(j)= \beta_{L} \mu(j),
\end{eqnarray*}
where each number is an index of $V_{0}$ thus implies each vertex. Recall that $\mu$ was the Perron Frobenius eigenvector of the original connection, and that $\om=[[M:N]]$. 
\end{proposition}
The above proposition implies that $[Q:P]= \om^{s-1}$.

 The following lemma is not necessary but noteworthy. For simplicity we omit $\otimes_{N}$ as long as there is no confusion.  
\begin{lemma}  
For the set of $N$-$N$ bimodules ${\cal X}:=\{X_{k}\}$ and any $Y \in \X$, $s \in \N$,  the following equality holds: 
$$ \sum_{X_{i} \in \X } N_{X_{1}, ... , X_{s}}^{Y} \mu_{1}...\mu_{s}=\om^{s-1} \mu_{Y},$$
where $N_{X_{1}, ... , X_{n}}^{Y}:={\rm dim Hom}(X_{1} ... X_{n}, Y)$, and $\mu_{i}=\mu(X_{i})$. 
\end{lemma}
{\bf Proof.}\\
We proceed by induction. The case $s=2$ is shown in Lemma 12.10 in \cite{EK7}. Suppose it holds for $s-1$. Note that $N_{X_{1}, ... , X_{s}}^{Y} = \sum_{Z} 
N_{X_{1}, ... , X_{s-1}}^{Z} N_{Z, X_{s}}^{Y}.$ Thus 
\begin{eqnarray*}
&& \sum_{X_{i} \in \X } N_{X_{1}, ... , X_{s}}^{Y} \mu_{1}...\mu_{s} \\
&=& \sum_{X_{s}, Z \in \X} N_{Z, X_{s}}^{Y} \mu_{s} 
\sum_{X_{i} \in \X } N_{X_{1}, ... , X_{s-1}}^{Z} \mu_{1}...\mu_{s-1} \\
&=& \sum_{X_{s}, Z \in \X} N_{Z, X_{s}}^{Y} \mu_{s}  \mu(Z) \om^{s-2}  \ \ \mbox{\rm (by the inductive hypothesis)} \\
&=& \om^{s-1} \mu(Y) \ \ \mbox{\rm (using the case} \ s=2)
\end{eqnarray*}
 QED. 

{\it Department of Mathematics,} \\
{\it University of California, Riverside, } \\
{\it 900 Big Springs Dr.} \\
{\it Riverside, CA 92521, USA,} \\
{\it email: marta@math.ucr.edu}

\end{document}